\documentclass[12pt, reqno]{amsart}

\usepackage{amsmath}
\usepackage{amssymb,amsfonts,amscd,amsthm}
\usepackage{latexsym}
\usepackage{mathrsfs}
\usepackage{amsaddr}
\usepackage{color}

\usepackage{hyperref}

\newcommand\de{\delta}

\newcommand\Om{\Omega}

\newcommand\R{\mathbb R}
\newcommand\C{\mathbb C}

\newcommand\Z{\mathbb Z}

\newcommand\gl{\mathfrak{gl}}

\newcommand\gr{\operatorname{gr}}

\newcommand\Mat{\operatorname{Mat}}

\newcommand\Shift{\operatorname{Shift}}
\newcommand\Comp{\operatorname{Comp}}

\newcommand\U{\mathcal U}

\newcommand\Al{\mathrm R}
\newcommand\M{\mathrm M}

\newcommand\AL{\mathcal  A}
\newcommand\NA{\mathcal N_{\AL}}
\newcommand\T{\mathrm T}
\newcommand\Tc{\T_{\textrm cycl}}
\newcommand\bl{\{\!\!\{}
\newcommand\br{\}\!\!\}}
\newcommand\x{\mathbf{x}}
\newcommand\y{\mathbf{y}}

\newcommand\aaa{\mathbf a}
\newcommand\bbb{\mathbf b}

\newcommand\coag{\rhd}

\newtheorem{theorem}{Theorem}[section]
\newtheorem{proposition}[theorem]{Proposition}
\newtheorem{lemma}[theorem]{Lemma}
\newtheorem{lemma-def}[theorem]{Lemma-definition}
\newtheorem{corollary}[theorem]{Corollary}

\theoremstyle{definition}
\newtheorem{definition}[theorem]{Definition}
\newtheorem{remark}[theorem]{Remark}

\numberwithin{equation}{section}

\begin{document}

\title[]{Remarks on Yangian-type algebras and double Poisson brackets}

\author{Grigori Olshanski and Nikita Safonkin}

\thanks{Supported by the Russian Science Foundation under project 23-11-00150}

\begin{abstract}
This short note is an announcement of results. We continue the study of Yangian-type algebras initiated in \cite{Ols}. These algebras share a number of properties of the Yangians of type A but are more massive. We refine and substantially enlarge the construction of \cite{Ols}. A direct link with the class of  linear double Poisson brackets on the free associative algebras (Pichereau and Van de Weyer, 2008) is also established.
\end{abstract}

\date{}

\maketitle

\tableofcontents

\section{Introduction}\label{sect1}

The present note is a continuation of \cite{Ols}. In that paper, a two-parameter family $\{Y_{d,L}\colon d,L=1,2,\dots\}$  of algebras was introduced and studied. For the particular value $L=1$ these algebras are the well-known Yangians $Y(\gl(d,\C))$ of the general Lie algebras $\gl(d,\C)$.  For $L\ge2$, the algebras $Y_{d,L}$ are new objects.  

It is well known that the Yangian $Y(\gl(d,\C))$ is a deformation of $\U(\gl(d,\C[u]))$, the universal enveloping algebra of the Lie algebra $\gl(d,\C[u])=\gl(d,\C)\otimes\C[u]$; the latter is a polynomial current Lie algebra. 

As shown in \cite{Ols}, a similar fact holds for $Y_{d,L}$: it is a deformation of $\U(\gl(d,\Al_L))$, where $\{\Al_L\colon L=1,2,3,\dots\}$ is a certain family of $\Z_{\ge0}$-graded associative algebras. However, in contrast with the Yangian case, the algebra $\Al_L$ with $L\ge2$ is non-commutative. Further, it is much more massive than $\C[u]$; namely, its $n$th homogeneous component ($n=0,1,2,\dots$) has dimension $d^2 L^{n+1}$ --- a quantity that grows exponentially in $n$. Nevertheless, the algebras $Y_{d,L}$ with $L\ge2$ share a number of properties of $Y(\gl(d,\C))$. For this reason, these new algebras were called \emph{Yangian-type algebras}. 

The purpose of the present note is to show that the results of \cite{Ols} in fact hold in a greater generality --- for a family  $\{Y_d(\Om)\}$ of  Yangian-type algebras, which depend (besides the previous parameter $d=1,2,3,\dots$) on an arbitrary  associative $\C$-algebra $\Om$. In this picture, the algebras $Y_{d,L}$ correspond to the particular case of $\Om=\C^{\oplus L}$, the direct sum of $L$ copies of the field $\C$.  

We also establish a connection between the algebras $Y_d(\Om)$ and a class of double Poisson brackets (in the sense of Van den Bergh \cite{vdB}), introduced in the work \cite{PW} by Pichereau and van de Weyer. 

Note that the correspondence $\Om\rightsquigarrow Y_d(\Om)$ is functorial in the sense that any algebra morphism $\Om\to\Om'$ entails an algebra morphism $Y_d(\Om)\to Y_d(\Om')$, for each fixed $d$. On the other hand, for fixed $\Om$ and  varying $d$, the algebras form a nested chain:
$$
Y_1(\Om)\subset Y_2(\Om)\subset Y_3(\Om)\subset\dots\,.
$$

As in \cite{Ols}, the algebra $Y_d(\Om)$ is extracted from a larger algebra $A_d(\Om)$,  which in turn is obtained through a projective limit construction: $A_d(\Om)=\varprojlim A_d(\Om;N)$, where $A_d(\Om;N)$ are certain centralizers in the universal enveloping algebras $\U(\gl(N,\Om))$. This is a version of the so-called centralizer construction (for the old version, related to the Yangians, see \cite[sect. 2.1]{Ols-Limits}, \cite[ch. 8]{M}). 

\medskip

The structure of the note is the following.

\smallskip

$\bullet$ Section \ref{sect2} introduces the centralizer construction leading to the algebras $A_d(\Om)$. 

\smallskip

$\bullet$ In section \ref{sect3} we construct some special elements of $A_d(\Om)$ (Theorem \ref{thm3.A}).

\smallskip

$\bullet$ In section \ref{sect4}, we define the algebra $Y_d(\Om)$ as the subalgebra of $A_d(\Om)$ generated by these special elements.  Then we formulate four theorems, which we consider as our main results.

--- Theorem \ref{thm4.A} states that, as a vector space, $A_d(\Om)$ splits into the tensor product of the subalgebras $A_0(\Om)$ and $Y_d(\Om)$; this shows to what extent $Y_d(\Om)$ differs from $A_d(\Om)$. Note that in the Yangian case, when $\Om=\C$,  these two subalgebras commute and the former subalgebra, $A_0(\C)$ is commutative. But in the general case, this is no longer true. 

--- Theorem \ref{thm4.B} is an analogue of the Poincar\'e-Birkhoff-Witt for $Y_d(\Om)$; it specifies the ``size'' of the algebra. 

--- Theorem \ref{thm4.C} provides a one-parameter group of automorphisms of $Y_d(\Om)$. In the Yangian case, these automorphisms are a useful tool for constructing representations. 

--- Theorem \ref{thm4.D} describes the structure of the graded algebra $\gr A_0(\Om)$ associated with $A_0(\Om)$. 

\smallskip

$\bullet$ In section \ref{sect5} we deal with Poisson brackets. By the very definition of the algebra $A_d(\Om)$, it has a canonical filtration such that associated graded algebra $\gr A_d(\Om)$ is commutative and carries a Poisson bracket of degree $-1$. Consequently, the algebras  $\gr Y_d(\Om)$ and $\gr A_0(\Om)$ are also graded Poisson algebras. Their structure is described in Propositions \ref{prop5.C} and \ref{prop5.D}, and it turns out that it fits into the general formalism of Van den Bergh \cite{vdB}. 

In more detail, recall that a double Poisson bracket $\bl-,-\br$ on an associative algebra $\AL$ is a bilinear map $\AL\times\AL\to \AL\otimes\AL$ subject to three axioms  \cite{vdB} (some versions of the skew-symmetry condition, the Leibniz rule, and the Jacobi identity). The paper  \cite{vdB} contains (among many other things) two general constructions. First, $\bl-,-\br$ determines a Poisson structure on the space of finite-dimensional representations of the algebra $\AL$. Second,  $\bl-,-\br$ determines a Lie algebra structure on the vector space $\AL/[\AL,\AL]$.  

We observe that the Poisson brackets computed in Propositions \ref{prop5.C} and \ref{prop5.D} are linked with these two constructions, and the corresponding double Poisson brackets  (we denote them by $\bl-,-\br_\Om$) are precisely the  ``linear'' double brackets on free associative algebras; a classification of such kind of brackets was given by Pichereau and Van de Weyer  \cite{PW}. As shown in \cite{PW}, they are specified by an arbitrary collection of structure constants of an associative algebra; in our context this is the algebra  $\Om$.

\smallskip

$\bullet$ In section \ref{sect6} we assign to the algebra $\Om$ a graded vector space, 
\begin{equation}\label{eq1.B}
\Al(\Om)=\bigoplus_{m=0}^\infty \Al_m(\Om), \qquad \Al_m(\Om):=\Om^{\otimes (m+1)},
\end{equation} 
and endow it with an associative multiplication 
\begin{equation}\label{eq1.A}
\Al_m(\Om)\times \Al_n(\Om)\to \Al_{m+n}(\Om), \quad m,\, n\in\Z_{\ge0},
\end{equation}
making $\Al(\Om)$ a graded associative algebra. Let us emphasize that the operation \eqref{eq1.A} differs from the standard product in the tensor algebra $\T(\Om)$. For the standard multiplication of tensors on the space $\Om$ we do not need the algebra structure on  $\Om$, whereas the operation \eqref{eq1.A} uses this structure. It is also important that the grading in \eqref{eq1.B} is shifted by $1$ with respect to the standard grading in $\T(\Om)$. A lucid interpretation of the algebra $\Al(\Om)$ is given in Proposition \ref{rem6.B}.

In the particular case of $\Om=\C^{\oplus L}$, one obtains the algebra $\Al_L$ mentioned in the beginning of the introduction (note that it is isomorphic to the path algebra of a quiver). 

We show (Proposition \ref{prop6.A}) that $Y_d(\Om)$ is a deformation of the universal enveloping algebra $\U(\gl(d,\Al(\Om))$. This extends the result from \cite{Ols} (mentioned above) that $Y_{d,L}$ is a deformation of $\U(\gl(d,\Al_L))$. Thus, the algebras  $\Al(\Om)$ may be viewed as a noncommutative analogue of polynomials in the construction of current algebras.

\section{The centralizer construction}\label{sect2}

Throughout the whole paper $\Om$ stands for an arbitrary associative algebra  of finite or countable dimension over $\C$ (possibly not unital). 

\subsection{The Lie algebra $\gl(N,\Om)$}
Let $\Mat(N,\C)$ denote the associative algebra of $N\times N$ matrices over $\C$. Likewise, $\Mat(N,\Om)$ is the algebra of $N\times N$ matrices with entries in $\Om$. Equivalently,
$$
\Mat(N,\Om):=\Mat(N,\C)\otimes \Om.
$$
The Lie algebra $\gl(N,\Om)$ coincides with $\Mat(N,\Om)$ as a vector space, and the bracket is given by the commutator. 

We denote by $E_{ij}$ the matrix units in $\Mat(N,\C)$  ($1\le i,j\le N$). They form a basis of $\gl(N,\C)$. For $x\in\Om$ we set 
$$
E_{ij}(x):=E_{ij}\otimes x \in \gl(N,\Om).
$$

The following relations hold
$$
[E_{ij}(x),E_{kl}(y)]=\de_{kj}E_{il}(xy)-\de_{il} E_{kj}(yx), \quad x,y\in\Om, \quad i,j,k,l\in\{1,\dots,N\}.
$$

The elements $E_{ij}(x)$ with $i,j\in\{1,\dots,N\}$ and $x$ ranging over a given basis of $\Om$ form a basis in $\gl(N,\Om)$.

\subsection{The action of $\gl(N,\C)$ on $\gl(N,\Om)$}

 If $\Om$ has a unit $1$, then $\gl(N,\C)$ is embedded into $\gl(N,\Om)$ via the map $E_{ij}\mapsto E_{ij}(1)$. Then we can restrict the adjoint action of $\gl(N,\Om)$ to the subalgebra $\gl(N,\C)$ and obtain in this way an action of $\gl(N,\C)$ by derivations of the Lie algebra $\gl(N,\Om)$.
 
 However, this action is well defined even if $\Om$ does not have a unit. Indeed, we define the action as $\operatorname{ad}_{\gl(N,\C)}\otimes\operatorname{Id}$, where $\operatorname{ad}_{\gl(N,\C)}$ denotes the adjoint representation of the Lie algebra $\gl(N,\C)$ and we identify $\gl(N,\Om)$ with $\gl(N,\C)\otimes\Om$.  

We agree to denote the action of $E_{ij}$ as $[E_{ij}, -]$, even if $\Om$ does not contain a unit. In this notation, we have
\begin{equation*}
[E_{ij}, E_{kl}(y)]=\de_{kj} E_{il}(y)- \de_{il}E_{kj}(y), \quad y\in\Om.
\end{equation*}

\subsection{The subalgebras $\gl_d(N,\C)$ and centralizers}

For $d=0,1,\dots,N$, we set
$$
\gl_d(N,\C):=\text{linear span of the matrix units $E_{ij}$ with $d+1\le i,j\le N$}.
$$
This is a Lie subalgebra of $\gl(N,\C)$, isomorphic to $\gl(N-d,\C)$. 

Consider $\U(\gl(N,\Om))$, the universal enveloping algebra of the Lie algebra $\gl(N,\Om)$. We set
\begin{equation*}
A_d(\Om;N):= \U(\gl(N,\Om))^{\gl_d(N,\C)}.
\end{equation*}
That is, $A_d(\Om;N)\subset \U(\gl(N,\Om))$ is the subalgebra of invariants of the action of the  Lie subalgebra $\gl_d(N,\C)\subset\gl(N,\C)$. 

When $\Om$ is unital, $\gl_d(N,\C)$ becomes a subalgebra of $\gl(N,\Om)$ and then $A_d(\Om;N)$ is the centralizer of this subalgebra. With a slight abuse of terminology, we call $A_d(\Om;N)$ a centralizer even when $\Om$ is non-unital. 

\subsection{Canonical projections $\pi_{N,N-1}$}

We denote by $I^+(N)$ the left ideal in $\U(\gl(N,\Om))$ generated by the elements $E_{iN}(x)$, where $i=1,\dots,N$ and $x\in\Om$. Likewise, $I^-(N)$ is the right ideal in $\U(\gl(N,\Om))$ generated by the elements $E_{Nj}(x)$, where $j=1,\dots,N$ and $x\in\Om$. 

\begin{lemma}
Let $\U(\gl(N,\Om))^{E_{NN}}=\U(\gl(N,\Om))^{\gl_{N-1}(N,\C)}$ be the centralizer of $E_{NN}$: that is, the subalgebra of elements annihilated by $[E_{NN},-]$. 

{\rm(i)} We have
\begin{equation*}
\U(\gl(N,\Om))^{E_{NN}}\cap I^+(N)=\U(\gl(N,\Om))^{E_{NN}}\cap I^-(N).
\end{equation*}

{\rm(ii)} Denote this intersection by $L(N)$. It is a two-sided ideal in $\U(\gl(N,\Om))^{E_{NN}}$, and the following direct sum decomposition holds
\begin{equation}\label{eq2.B}
\U(\gl(N,\Om))^{E_{NN}}=\U(\gl(N-1,\Om))\oplus L(N).
\end{equation}
\end{lemma}

Let 
$$
\pi_{N,N-1}: \U(\gl(N,\Om))^{E_{NN}}\to \U(\gl(N-1,\Om))
$$
be the projection onto the first summand in \eqref{eq2.B}. Since $L(N)$ is a two-sided ideal, $\pi_{N,N-1}$ is an algebra homomorphism. It is consistent with the canonical filtration in universal enveloping algebras. 

\subsection{The projective limit algebras $A_d(\Om)$}

\begin{lemma}
Let $d=0,1,2,\dots$ and $N>d$. Then
$$
\pi_{N,N-1}(A_d(\Om;N)) \subseteq A_d(\Om;N-1).
$$
\end{lemma}

This lemma makes it possible to give the following definition.

\begin{definition}[The centralizer construction]
For $d=0,1,2,\dots$, we set
\begin{equation*}
A_d(\Om):=\varprojlim (A_d(\Om;N), \pi_{N,N-1}), 
\end{equation*}
where the projective limit is taken in the category of filtered algebras. Thus, the limit algebra $A_d(\Om)$ is a filtered algebra, too. 
\end{definition}

By the very definition, every element $\mathbf a\in A_d(\Om)$ is a sequence $\{\mathbf a(N): N\ge d\}$ such that
$$
\text{$\mathbf a(N)\in A_d(\Om;N)$ for $N\ge d$, \quad $\pi_{N,N-1}(\mathbf a(N))=\mathbf a(N-1)$ for $N>d$},
$$
and $\deg \mathbf a(N)$ is bounded by a constant independent  of $N$. 

We denote by $\pi_{\infty,N}$ the natural projection $A_d(\Om)\to A_d(\Om;N)$. In this notation, $\pi_{\infty,N}(\aaa)=\aaa(N)$. 

There are natural embeddings
\begin{equation}\label{eq2.A}
A_0(\Om)\subset A_1(\Om)\subset A_2(\Om)\subset\dots\,.
\end{equation}
In particular $A_0(\Om)$ is a subalgebra of $A_d(\Om)$ for any $d\in\Z_{\ge1}$.

\section{Special elements of $A_d(\Om)$}\label{sect3}

\subsection{The elements $e_{ij}(\x;N)$}

For $m=1,2,\dots$, let $\Om^m$ be the $m$-fold direct product $\underbrace{\Om\times\dots\times\Om}_m$. Thus, an element of $\Om^m$ is an ordered $m$-tuple $\x=(x_1,\dots,x_m)$ of elements from $\Om$. 

Given an element $\x=(x_1,\dots,x_m)\in\Om^m$, a pair of indices $i,j$, and an integer $N\ge\max(i,j)$, we define an element of $\U(\gl(N,\Om))$ as follows:
\begin{equation}
e_{ij}(\x;N):=\sum_{a_1,\dots,a_{m-1}=1}^N E_{ia_1}(x_1) E_{a_1a_2}(x_2)\dots E_{a_{m-1}j}(x_m)
\end{equation}
with the understanding that 
$$
e_{ij}(\x;N)=E_{ij}(x_1)\quad \text{when $m=1$ and $\x=(x_1)$.}
$$

\subsection{The coagulation transforms $\x\mapsto\x\coag\nu$}
Recall that a \emph{composition} is an ordered tuple $\nu=(\nu_1,\dots,\nu_n)$ of positive integers. We set 
$$
|\nu|:=\nu_1+\dots+\nu_n, \quad \ell(\nu)=n,
$$
and denote by $\Comp(m)$ the set of compositions $\nu$ with $|\nu|=m$. 

\begin{definition}
To an $m$-tuple  $\x=(x_1,\dots,x_m)\in\Om^m$ and a composition $\nu\in\Comp(m)$ with $\ell(\nu)=n$, we assign an $n$-tuple $y:=\x\coag\nu=(y_1,\dots,y_n)\in\Om^n$ as follows:
$$
y_1=x_1\dots x_{\nu_1}, \quad y_2=x_{\nu_1+1}\dots x_{\nu_1+\nu_2}, \quad\dots\quad y_n=x_{\nu_1+\dots+\nu_{n-1}+1}\dots x_m,
$$
where the products are taken in the algebra $\Om$. Let us call this transform \emph{coagulation}. 
\end{definition}

For instance, for $m=3$ there are 4 compositions,
$$
\nu=\; (1,1,1), \; (2,1), \; (1,2), \; (3),
$$
and the corresponding elements are
$$
\y=\; (x_1,x_2,z_3), \; (x_1x_2, x_3), \; (x_1, x_2x_3), \; x_1x_2x_3.
$$

\subsection{The elements $t_{ij}(\x;N;s)$ and connection coefficients}

Set
\begin{equation}\label{eq3.A}
W(\Om):=\bigsqcup_{m=1}^\infty \Om^m.
\end{equation}
For $\x\in\Om^m\subset W(\Om)$ we set $\ell(\x):=m$. 

\begin{definition}
Let $s\in\C$ be an auxiliary parameter. As before, let $i,j,N$ be positive integers, $1\le i,j\le N$.  Next, let $\x\in W(\Om)$ and $m:=\ell(\x)$. From these data we define elements of the algebra $\U(\gl(N,\Om))$, as follow
\begin{equation*}
t_{ij}(\x;N;s):=\sum_{\nu\in\Comp(m)} (-N-s)^{\ell(\x)-\ell(\x\coag\nu)}e_{ij}(\x\coag\nu;N).
\end{equation*}
\end{definition}

Here we agree that $0^0:=1$. From this convention it follows that 
$$
t_{ij}(\x;N;s)\big|_{s=-N}=e_{ij}(\x;N).
$$

Note that $\x\coag\nu=\x$ for the composition $\nu=(1,\dots,1)$, while for all other compositions, $\ell(\x\coag\nu)<\ell(x)$. It follows that 
$$
t_{ij}(\x;N;s)=e_{ij}(\x;N)+(\cdots),
$$ 
where the dots stand for a linear combination of elements of the form $e_{ij}(\y;N)$ with $\ell(\y)<\ell(\x)=m$.

The next lemma connects elements corresponding to different values of the parameter $s$. 

\begin{lemma}\label{lemma3.A}
For any $s,s'\in\C$,
\begin{equation*}
t_{ij}(\x;N;s)=\sum_{\nu\in\Comp(m)}(s'-s)^{\ell(\x)-\ell(\x\coag\nu)} t_{ij}(\x\coag\nu;N;s').
\end{equation*}
\end{lemma}

\subsection{The elements $t_{ij}(\x;s)\in A_d(\Om)$}

\begin{theorem}\label{thm3.A}
Fix $\x\in W(\Om)$ and $s\in\C$. One has
$$
\pi_{N,N-1}(t_{ij}(\x;N;s))=t_{ij}(\x;N-1;s), \quad N>\max(i,j).
$$
\end{theorem}

\begin{corollary}\label{cor3.A}
Fix $d\in\Z_{\ge1}$. For any $\x\in W(\Om)$, $s\in\C$, and $i,j\in\{1,\dots,d\}$ there exists an element $t_{ij}(\x;s)\in A_d(\Om)$, uniquely determined by the property that 
$$
\pi_{\infty,N}(t_{ij}(\x;s))=t_{ij}(\x;N;s), \qquad \forall N\ge d.
$$
\end{corollary}

Together with Lemma \ref{lemma3.A} this in turn implies

\begin{corollary}\label{cor3.B}
For any $s,s'\in\C$,
\begin{equation*}
t_{ij}(\x;s)=\sum_{\nu\in\Comp(m)}(s'-s)^{\ell(\x)-\ell(\x\coag\nu)} t_{ij}(\x\coag\nu;s').
\end{equation*}
\end{corollary}

\section{Main results}\label{sect4}

\begin{definition}
Given $d\in\Z_{\ge1}$, we denote by $Y_d(\Om)$ the unital subalgebra of $A_d(\Om)$ generated by the elements of the form $t_{ij}(\x;s)$, where $i,j\in\{1,\dots,d\}$, $\x\in W(\Om)$, and $s\in\C$ is an arbitrary fixed number. By Corollary \ref{cor3.B}, the definition does not depend on the choice of $s$.
\end{definition}

Fix a basis $X$ in $\Om$ and set (cf. \eqref{eq3.A})
\begin{equation*}
W(X):=\bigsqcup_{m=1}^\infty X^m \subset W(\Om).
\end{equation*}

Because an element $t_{ij}(\x; s)$ with $\x=(x_1,\dots,x_m)$ depends multilinearly on $x_1,\dots,x_m$, it is a linear combination of some elements $t_{ij}(\y;s)$ with $\y\in X^m$. It follows that in the definition above we may assume that $\x$ ranges over $W(X)$. 

Pick an arbitrary total order $\prec$ on the set of triples $(i,j,\x)$, where $i,j\in\{1,\dots,d\}$ and $\x\in W(X)$. We say that a product
\begin{equation}\label{eq4.A}
t_{i_1j_1}(\x^1; s)\dots t_{i_rj_r}(\x^r;s), \quad \text{where $\x^1,\dots,\x^r\in W(X)$}, 
\end{equation}
is an \emph{ordered monomial} if the corresponding triples $(i_1,j_1,\x^1), \dots, (i_r,j_r,\x^r)$ weakly increase with respect to $\prec$. 

\begin{theorem}[Splitting of $A_d(\Om)$]\label{thm4.A}
Let $m\colon A_d(\Om)\times A_d(\Om)\to A_d(\Om)$ denote the multiplication map. The subalgebras $A_0(\Om)$ and $Y_d(\Om)$ split $A_d(\Om)$ in the sense that 
$$
m\colon A_0(\Om)\otimes Y_d(\Om)\to A_d(\Om)
$$
is an isomorphism of vector spaces.
\end{theorem}

\begin{theorem}[Analogue of the Poincar\'e-Birkhoff-Witt theorem for $Y_d(\Om)$]\label{thm4.B}
For any choice of a basis $X\subset \Om$, a total order $\prec$ as above, and any fixed $s\in\C$, the ordered monomials \eqref{eq4.A} together with $1$ form a basis of the algebra $Y_d(\Om)$.
\end{theorem}

\begin{theorem}[Shift automorphisms of the algebra $Y_d(\Om)$]\label{thm4.C}
The additive group $\C$ acts on the algebra $Y_d(\Om)$ by automorphisms $\Shift_c$, $c\in\C$, such that 
$$
\Shift_c\colon t_{ij}(\x;s)\mapsto t_{ij}(\x;s+c)
$$
for any $i,j\in\{1,\dots,d\}$, $\x\in W(\Om)$, and $s\in\C$. 
\end{theorem}

Introduce a notation: 
\begin{equation*}
\T^+(\Om):=\bigoplus_{m=1}^\infty \T^m(\Om),  \quad \T^m(\Om):=\Om^{\otimes m},
\end{equation*}
that is, $\T^+(\Om)$ is the space of tensors over $\Om$ of strictly positive degree;
\begin{equation*}
\M_d(\Om):=\Mat(d,\C)\otimes \T^+(\Om);
\end{equation*}
$\mathcal S(\M_d(\Om))$ is the symmetric algebra over the vector space $\M_d(\Om)$. We equip $\M_d(\Om)$ with the grading induced by the grading of the tensor space  $\T^+(\Om)$, and next extend the grading from $\M_d(\Om)$ to the symmetric algebra $\mathcal S(\M_d(\Om))$.  

Theorem \ref{thm4.B} implies

\begin{corollary} 
There is an isomorphism of graded algebras 
\begin{equation}\label{eq4.B}
\gr Y_d(\Om)\simeq \mathcal S(\M_d(\Om)).
\end{equation}
\end{corollary} 

Now denote by $\Tc^m(\Om)$ the space of coinvariants for the action of the group $\Z/m\Z$ by cyclic permutations on the space $\T^m(\Om)$ ($m=1,2,\dots$). Next, introduce the graded vector space 
$$
\Tc^+(\Om):=\bigoplus_{m=1}^\infty \Tc^m(\Om)
$$
and denote by $\mathcal S(\Tc^+(\Om))$ the symmetric algebra over this space. 

\begin{theorem}\label{thm4.D}
There is an isomorphism of graded algebras  
\begin{equation}\label{eq4.C}
\gr A_0(\Om)\simeq \mathcal S(\Tc^+(\Om)).
\end{equation}
\end{theorem}

\section{Link with double Poisson brackets}\label{sect5}

Due to commutativity of the algebras $\gr Y_d(\Om)$ and $\gr A_0(\Om)$ they possess canonical Poisson brackets. We will provide an explicit description of these brackets using the isomorphisms \eqref{eq4.B} and \eqref{eq4.C}. This will reveal a connection with double Posson brackets in the sense of Van den Bergh \cite{vdB}.

\subsection{Double Poisson brackets \cite{vdB}}

Let $\AL$ be an arbitrary associative algebra. We need its tensor square, which will be denoted by the symbol  $\AL^{\boxtimes 2}$. Here we do not use the standard symbol $\otimes$ for the tensor multiplication, because a bit later the algebra $\AL$ will itself consist of tensors.  

Recall that a double Poisson bracket $\bl-,-\br$ on $\AL$ is a bilinear map 
$$
\AL\times\AL \to \AL^{\boxtimes2}, \quad (a,b) \mapsto \bl a,b\br,
$$
subject to three axioms, which are analogues of the skew-symmetry condition, the Leibniz rule, and the Jacobi identity (Van den Bergh \cite{vdB}). 

Note that if $\AL$ contains a unit $1$, then $\bl 1,a\br=\bl a,1\br=0$ for all $a\in\AL$.

\subsection{Linear double Poisson brackets \cite{PW}}

Assume for a moment that $\Om$ is simply a vector space without any additional structure. Denote by $\T(\Om)$ the tensor algebra over $\Om$, or, which is the same, the free algebra generated by the vector space $\Om$. From now on we set $\AL:=\T(\Om)$.  

Let $\bl-,-\br$ be a double poisson bracket on $\T(\Om)$. By one of the axioms (analogue of the Leibniz rule), $\bl-,-\br$  is uniquely determined by its values on $\Om\times\Om$, where we identify $\Om$ with $\T^1(\Om)$. Following \cite{PW}, we say that $\bl-,-\br$ is  \emph{linear} if its restriction to  $\Om\times\Om$ takes values in  $1\boxtimes \Om +\Om\boxtimes 1$. Then we necessarily obtain that  
\begin{equation}\label{eq5.A}
\bl x,y\br = 1\boxtimes \mu(x,y) - \mu(y,x)\boxtimes 1, \quad x,y\in\Om,
\end{equation} 
for a certain bilinear map $\mu: \Om\times\Om\to \Om$. However, $\mu$ cannot be arbitrary (because an analogue of the Jacobi identity must hold).

\begin{proposition}[Pichereau and Van de Weyer \cite{PW}, Proposition 10]\label{prop5.A}
In order for the formula \eqref{eq5.A} to determine a double Poisson bracket it is necessary and sufficient that $\mu$ be an associative multiplication on the vector space $\Om$.
\end{proposition}
(See also Odesskii, Rubtsov, and Sokolov \cite[(2.9)]{ORS}.)
Thus, for an arbitrary vector space $\Om$ there is a one-to-one correspondence between the linear double Poisson brackets on the free algebra $\T(\Om)$ and the associative algebra structures on $\Om$. (In \cite{PW} and \cite{ORS}, the authors deal with finite-dimensional vector spaces, but in  Proposition \eqref{prop5.A}, this restriction is unessential.)

\subsection{The double Poisson bracket $\bl-,-\br_\Om$ on $\T(\Om)$}
 
We assume again that $\Om$ is an associative algebra. Let $\bl-,-\br_\Om$ denote the double Poisson bracket on $\T(\Om)$ afforded by Proposition \ref{prop5.A}. Here is an explicit expression for this bracket.  

Introduce first a new notation. 
Let $\x=x_1\otimes\dots\otimes x_m\in\T^m(\Om)$ and $\y=y_1\otimes\dots\otimes y_n\in\T^n(\Om)$ be two non-scalar decomposable tensors. For $r=1,\dots,m$ and $s=1,\dots,n$, we set  
\begin{align*}
\x^{<r}&=x_1\otimes \dots \otimes x_{r-1},  \quad \x^{>r}=x_{r+1}\otimes \dots\otimes x_m, \\
\y^{<s}&=y_1\otimes\dots \otimes y_{s-1},  \quad \y^{>s}=y_{s+1}\otimes \dots\otimes y_n,
\end{align*}
with the undestanding that $\x^{<1}$, $\x^{>m}$, $\y^{<1}$, and $\y^{>n}$ are interpreted as the unit $\mathbf1\in\T(\Om)$.  

\begin{proposition}\label{prop5.B}
The double Poisson bracket $\bl-,-\br_\Om$ on the free algebra $\T(\Om)$, defined according to the general recipe of Proposition \ref{prop5.A} is given by the formula 
\begin{equation}\label{eq5.B}
\begin{aligned}
\bl\x, \y\br_\Om &=\sum_{r=1}^m\sum_{s=1}^n 
\bigg( (\y^{<s}\otimes\x^{>r})\boxtimes (\x^{<r}\otimes (x_r y_s)\otimes \y^{>s})\bigg) \\
&-\sum_{r=1}^m\sum_{s=1}^n 
\bigg( (\y^{<s}\otimes (y_s x_r)\otimes \x^{>r})\boxtimes (\x^{<r}\otimes\y^{>s})\bigg).
\end{aligned}
\end{equation}
Finally, $\bl-,-\br_\Om$ is equal to $0$ if at least one of the arguments is in $\T^0(\Om)=\C$. 
\end{proposition}

\subsection{The Poisson bracket on $\gr Y_d(\Om)$}

Given $\x\in\T^+(\Om)$ and a bi-index $ij$, where $i,j\in\{1,\dots,d\}$, we set
$$
p_{ij}(\x):= E_{ij}\otimes\x\in \M_d(\Om).
$$

On the algebra $\mathcal S(\M_d(\Om))$ there exists a Poisson bracket $\{-,-\}_{d,\Om}$ uniquely determined by the property that  for decomposable tensors $\x\in\T^m(\Om)$, $\y\in\T^n(\Om)$, and indices $i,j,k,l\in\{1,\dots,d\}$, 
\begin{equation*}
\begin{aligned}
\{p_{ij}(\x), p_{kl}(\y)\}_{d,\Om}&=\sum_{r=1}^m\sum_{s=1}^n 
p_{kj} (\y^{<s}\otimes\x^{>r})p_{il} (\x^{<r}\otimes (x_r y_s)\otimes \y^{>s}) \\
&-\sum_{r=1}^m\sum_{s=1}^n 
p_{kj}(\y^{<s}\otimes (y_s x_r)\otimes \x^{>r})p_{il} (\x^{<r}\otimes\y^{>s}),
\end{aligned}
\end{equation*}
with the agreement that $p_{ab}(\mathbf1):=\de_{ab}$.

Writing the element $\bl\x,\y\br_\Om\in\T(\Om)^{\boxtimes2}$ in Sweedler notation as $\sum \bl\x,\y\br'\boxtimes \bl\x,\y\br''$, we can rewrite the formula as
\begin{equation}\label{eq5.C}
\{p_{ij}(\x), p_{kl}(\y)\}_{d,\Om}=\sum p_{kj}(\bl\x,\y\br')p_{il}(\bl\x,\y\br'').
\end{equation}

The transition from the double Poisson bracket $\bl-,-\br_\Om$ on the free algebra $\T(\Om)$ (formula \eqref{eq5.B}) to the Poisson bracket $\{-,-\}_{d,\Om}$ on the symmetric algebra $\mathcal S(\M_d(\Om))$ (formula \eqref{eq5.C}) ``almost'' fits into the Van den Bergh's formalism   \cite[Proposition 7.5.1]{vdB}. The reservation  ``almost'' is due to the fact that in our situation, we do not need to impose the relation  $\sum_q a_{pq}b_{qr}=(ab)_{pr}$, written at the bottom of the page 5751 in  \cite{vdB} (about this point, see  \cite[section 5.2]{Ols}).

\begin{proposition}\label{prop5.C}
Under the isomorphism \eqref{eq4.B} the canonical Poisson bracket on $\gr Y_d(\Om)$ coincides with the bracket $\{-,-\}_{d,\Om}$ given by the formula \eqref{eq5.C}. 
\end{proposition}

Thus, the algebra $Y_d(\Om)$ is a filtered deformation of the Poisson algebra 
$$(\mathcal S(\M_d(\Om)), \{-,-\}_{d,\Om}).$$ 

\subsection{The Poisson bracket on $\gr A_0(\Om)$}

Assume first that $\AL$ is an arbitrary associative algebra and let $[\AL,\AL]$ be its subspace spanned by the commutators $[\aaa,\bbb]$, where $\aaa,\bbb\in\AL$. Denote by $\NA$ the quotient space $\AL/[\AL,\AL]$. As shown by Van den Bergh \cite[Proposition 1.4 and Corollary 2.4.6]{vdB}, a double Poisson bracket $\bl-,-\br$ on $\AL$ gives rise to a Lie algebra structure on the space $\NA$ and hence to a Poisson algebra structure on the symmetric algebra $\mathcal S(\NA)$. 

Examine now the case when $\AL$ is the free algebra $\T(\Om)$ with the double Poisson bracket $\bl-,-\br_\Om$ (formula \eqref{eq5.B}). Then the space $\NA$ can be identified with $\C\oplus \Tc^+(\Om)$. Thus, $\C\oplus \Tc^+(\Om)$ acquires a Lie algebra structure. This Lie algebra is the direct sum of the trivial one-dimensional Lie algebra $\C$ and the Lie subalgebra $\Tc^+(\Om)$; we are interested in the latter subalgebra. The Lie algebra structure on $\Tc^+(\Om)$ determines a Poisson bracket on $\mathcal S(\Tc^+(\Om))$, let us denote it by $\{-,-\}_\Om$. 

\begin{proposition}\label{prop5.D}
Let, as usual, $\Om$ be an arbitrary associative algebra. Under the isomorphism \eqref{eq4.C}, the Poisson bracket on $\gr A_0(\Om)$ coincides with the Poisson bracket $\{-,-\}_\Om$ on $\mathcal S(\Tc^+(\Om))$ coming from the general Van den Bergh construction described above. 
\end{proposition}

Thus, the algebra $A_0(\Om)$ serves as a filtered quantization of the Poisson algebra 
$$
(\mathcal S(\Tc^+(\Om)), \{-,-\}_\Om).
$$
 
In the works Schedler \cite{Sch} and Ginzburg-Schedler  \cite{GinS}, the quantization problem for Poisson algebras of the form $\mathcal S(\NA)$ was studied in another context, related to quivers.

\section{Degeneration $Y_d(\Om)\rightsquigarrow \U(\gl(d,\Al(\Om)))$}\label{sect6}

As above, we work with an arbitrary associative algebra  $\Om$. Denote by $\Al(\Om)$ the vector space $\T^+(\Om)$ endowed with a new grading, starting from zero:
$$
\Al(\Om):=\bigoplus_{n=0}^\infty \Al(\Om)_n, \quad \Al(\Om)_n:=\T^{n+1}(\Om).
$$
Next, we endow $\Al(\Om)$ with a multiplication denoted by the symbol $\odot$: for decomposable tensors,  
$$
\x=x_1\otimes\dots\otimes x_{m+1}\in\Al_m(\Om), \quad \y=y_1\otimes\dots\otimes y_{n+1}\in\Al_n(\Om),  
$$
where  $m,n\in\{0,1,2,\dots\}$, the multiplication is as follows:
\begin{equation*}
\x\odot\y:=x_1\otimes\dots\otimes x_{m}\otimes (x_{m+1}y_1)\otimes  y_2\otimes \dots\otimes y_{n+1}\in\Al_{m+n}(\Om). 
\end{equation*}
this operation is associative and consistent with the new grading. Thus, $\Al(\Om)$ is a graded algebra.

\begin{remark} 

(i) The $0$th component of $\Al(\Om)$ is the initial algebra $\Om$. 

(ii) $\Al(\Om)$ has a unit if and only if this is so for $\Om$.

(iii)  $\Al(\Om)$ is non-commutative unless either $\Om=\C$ or the product in $\Om$ is null.

(iv) For $\Om=\C$, the corresponding algebra $\Al(\C)$ is isomorphic to the algebra of polynomials $\C[u]$: the isomorphism is given by the assignment  $1^{\otimes (r+1)}\mapsto u^r$.

(v) More generally, if $\Om=\C^{\oplus L}$ (the direct sum of $L$ copies of $\C$), then $\Al(\Om)$ is isomorphic to the path algebra of the quiver $Q_L$ with $L$ vertices indexed by $1,\dots,L$ and $L^2$ edges $i\to j$, one for each ordered pair of vertices $(i,j)$, including $i=j$.
\end{remark}

\begin{remark}\label{rem6.B}
The vector space $M:=\Om\otimes_{\mathbb{C}}\Om$ possesses a natural structure of a $\Om$-bimodule, so that we can built the graded tensor algebra over the bimodule $M$:
$$
\T_\Om(M):=\Om\oplus M \oplus (M\otimes_\Om M) \oplus (M \otimes_\Om M\otimes_\Om M) \oplus\dots,
$$
where all tensor products are taken over $\Om$. If $\Om$ contains a unit, then the algebra $\R(\Om)$ is isomorphic to the algebra $\T_\Om(M)$. We are grateful to Boris Feigin for this remark. It shows that the algebra $\R(\Om)$ is not something exotic. 
\end{remark}

Consider now the Lie algebra $\gl(d,\Al(\Om))$, where $d\in\{1,2,3,\dots\}$. As a vector space, it coincides with the associative algebra $\Mat(d,\C)\otimes\Al(\Om)$, and the Lie bracket is the commutator $[-,-]$ in this algebra. That is, 
$$
[X\otimes\x,Y\otimes\y]=(XY)\otimes (\x\odot\y)-(YX)\otimes(\y\odot\x), \quad X,Y\in\Mat(d,\C), \; \x,\y\in \Al(\Om).
$$
The grading of $\Al(\Om)$ induces, in a natural way, a grading in the Lie algebra $\gl(d,\Al(\Om))$, and from this we obtain a grading in the universal enveloping algebra $\U(\gl(d,\Al(\Om)))$. 

Along with the canonical filtration in $Y_d(\Om)$ we have dealt with so far, the algebra $Y_d(\Om)$ also admits another filtration, which we will call the ``shifted filtration": it is specified by the condition that  the degree of $t_{ij}(\x;s)$ equals not $\ell(\x)$ as before, but $\ell(\x)-1$. 

\begin{proposition}\label{prop6.A}
The graded algebra associated with the shifted filtration in $Y_d(\Om)$ is isomorphic to the algebra $\U(\gl(d,\Al(\Om)))$. 
\end{proposition}

\vspace{3ex}

\leftline{Grigori Olshanski}

\leftline{Institute for Information Transmission Problems, Moscow, Russia.} 
\leftline{Skolkovo Institute of Science and Technology, Moscow, Russia.} 
\leftline{HSE University, Moscow, Russia.}
\leftline{e-mail: olsh2007@gmail.com}

\vspace{3ex}

\leftline{Nikita Safonkin}
\leftline{Laboratoire de Math{\'e}matiques de Reims, }
\leftline{Universit{\'e} de Reims-Champagne-Ardenne, 
CNRS UMR 9008, F-51687, Reims, France.}
\leftline{Skolkovo Institute of Science and Technology, Moscow, Russia.} 
\leftline{e-mail: safonkin.nik@gmail.com}

\end{document}